\documentclass[12pt]{amsart}

\usepackage{amsmath}
\usepackage{url}
\usepackage{amsthm}
\usepackage{amssymb}
\usepackage{enumerate}
\usepackage{geometry}
\usepackage{hyperref}
\hypersetup{breaklinks=true,
pagecolor=white,
colorlinks=true,citecolor=black,filecolor=black,linkcolor=black,urlcolor=black
}

\newtheorem{thm}{Theorem}%[section]

\theoremstyle{definition}
\newtheorem{remark}[thm]{Remark}

\newcommand{\R}{\mathbb{R}}

\newcommand{\Z}{\mathbb{Z}}
\newcommand{\T}{\mathbb{T}}
\newcommand{\s}{\mathbb{S}}

\title{Character analysis using Fourier series}
\author{Daniel Reem}
\begin{document}
\maketitle
\vspace{-0.7cm}
\begin{abstract}
We show that character analysis using Fourier series is possible, 
at least when a mathematical character is considered. Previous approaches
to character analysis are somewhat not in the spirit of harmonic analysis.
\end{abstract}

\section*{The analysis}
Although it may seem somewhat ambitious, character analysis using Fourier  series is possible, at least when a mathematical character is considered, namely a homomorphism from a given group into the multiplicative group $\s^1$ of complex numbers whose modulus is one.

To understand how, consider the well-known and useful results in (abstract)  harmonic analysis saying that the character group (the group of continuous characters) of the (additive) flat torus $\T=\R/(2\pi\Z)$ is the group of integers $\Z$, and the  character group of the (additive) real line $\R$ is $\R$. Equivalently, any $2\pi$-periodic continuous function $f:\R\to\s^1$   satisfying $f(x+y)=f(x)f(y)$ for each $x,y\in \R$ must have the form $f(x)=e^{ik\cdot x}$ for some $k\in \Z$, and  any continuous homomorphism  $f:\R\to\s^1$ must have the form $f(x)=e^{i\alpha\cdot x}$ for some $\alpha\in \R$. 

Familiar proofs of these results \cite[pp. 79-80]{Deitmar}, \cite[pp. 143]{DeVito}, \cite[pp. 89-90]{Folland}, \cite[pp. 364-366]{HewittRoss}, \cite[pp. 44-46, 95-96]{Krantz}, \cite[pp. 246-247]{Pontryagin}, \cite[pp. 12-13]{RudinFourier}, \cite[pp. 231, 237]{SteinShakarchi} are usually restricted to the case of continuous mappings (although the arguments in \cite{Folland, Krantz, RudinFourier, SteinShakarchi} can be generalized), and, interestingly, they involve arguments which are somewhat not in the spirit of harmonic analysis. A natural tool in harmonic analysis is Fourier series. Using this tool, we can prove the following theorem: 
\begin{thm}\label{thm:PeriodicHomomorphism}
\begin{enumerate}[(a)]
\item Let $f:\R\to \s^1$ be $2\pi$-periodic. %, i.e., $f(x)=f(x+2\pi)$ for each $x\in \R$. 
Then $f$ is a  measurable homomorphism if and only if there exists $k\in \Z$ such that $f(x)=e^{ik\cdot x}$ for each $x\in \R$.
\item A function $f:\R\to \s^1$ is a  measurable homomorphism if and only if there exists $\alpha\in \R$ such that $f(x)=e^{i\alpha\cdot x}$ for each  $x\in \R$.
\end{enumerate}
\end{thm} 

\begin{proof}
\begin{enumerate}[(a)]
\item\label{item:Torus} It is evident that $f(x)=e^{ik\cdot x}$ is a $2\pi$-periodic measurable  homomorphism for each $k\in\Z$. On the other hand, let $f:\R\to \s^1$ be a $2\pi$-periodic measurable  homomorphism. Given any $k\in \Z$, the bounded function 
$f(x)e^{-i k\cdot x}$ is $2\pi$-periodic and integrable. Hence its integral on $I=[0,2\pi]$ is invariant under translations, i.e., $\int_I f(x+y)e^{-ik\cdot(x+y)}dx=\int_I f(x)e^{-ik\cdot x}dx$ for each $y\in \R$. Since $f$ is a homomorphism it satisfies $f(x+y)=f(x)f(y)$ for each $x,y\in \R$. Therefore its $k$th Fourier coefficient $\hat{f}(k)$ satisfies (for each $y$ in $\R$) 
\begin{equation*}
\hat{f}(k)=\frac{1}{2\pi}\int_I f(x)e^{-ik\cdot x} dx=\frac{1}{2\pi}\int_I f(x+y)e^{-ik\cdot (x+y)}dx
\\=\hat{f}(k)f(y)e^{-ik\cdot y}.
\end{equation*}
If $\hat{f}(k)=0$ for each $k\in \Z$, then $f(x)=0$ for almost every $x$ by the uniqueness of the Fourier expansion, contradicting the fact that $|f(x)|=1$ for each $x$. Thus $\hat{f}(k)\neq 0$ for some $k\in \Z$, and then $f(y)=e^{ik\cdot y}$ for all $y\in \R$, as claimed.

\item\label{item:R} It is evident that $f(x)=e^{i\alpha\cdot x}$ is a measurable  homomorphism for each $\alpha\in \R$. On the other hand, let $f:\R\to\s^1$  be a given measurable homomorphism. Since $|f(2\pi)|=1$ there exists $\beta\in [0,1)$ satisfying $e^{i 2\pi\beta}=f(2\pi)$. Define $g:\R\to \s^1$ by $g(x)=e^{i\beta\cdot x}$. Let $h=f/g$. Then $h$ is well defined and it is a quotient of measurable  homomorphisms. Hence $h$ is a measurable  homomorphism by itself. In addition, $h$ is $2\pi$-periodic,  since for each $x\in \R$
\begin{equation*}
h(x+2\pi)=\frac{f(x+2\pi)}{g(x+2\pi)}=\frac{f(x)f(2\pi)}{g(x) e^{i 2\pi \beta}}=\frac{f(x)}{g(x)}=h(x).
\end{equation*}
Consequently, by what we proved in part \eqref{item:Torus} it follows that $h(x)=e^{ik\cdot x}$ for some $k\in \Z$, and, as a result, $f(x)=e^{i\alpha\cdot x}$  for each $x\in \R$, where $\alpha=k+\beta$.
\end{enumerate}
\end{proof}
\begin{remark}
The above results and proofs can be generalized almost word for word to characters defined on the $n$-dimensional flat torus $\T^n$ and on $\R^n$. Now in the proof of Theorem \ref{thm:PeriodicHomomorphism}\eqref{item:Torus}  we consider periodic functions whose $n$ periods are $2\pi e_j$, for each  element $e_j,\,j=1,\ldots,n$, of the (standard) basis, and in the proof of Theorem  \ref{thm:PeriodicHomomorphism}\eqref{item:R} we take $\beta=(\beta_j)_{j=1}^n\in [0,1)^n$ satisfying $e^{i 2\pi \beta_j}=f(2\pi e_j)$ for each $j$.  
\end{remark}
\begin{remark}
The above results, although they  seem more general than the results described in the references, are in fact equivalent to them, since every measurable character defined on a locally compact group is continuous, as implied by 
\cite[Theorem ~22.18, p. 346]{HewittRoss}. However, the proof of this theorem is quite complicated, and in the cases considered in this note our proofs seem simpler. 
\end{remark}

{\bf \noindent Acknowledgment}\\
I thank the referee for a suggestion which helped to shorten this note.

\bigskip

\noindent\textit{Department of Mathematics, University of Haifa,  Haifa, Israel, and Department of Mathematics, The Technion - Israel Institute of Technology, Haifa, Israel. Current address: IMPA - Instituto Nacional de Matem\'atica Pura e Aplicada, Estrada Dona Castorina 110, Jardim Bot\^anico, CEP 22460-320, Rio de Janeiro, RJ,  Brazil. \\ 
Email: dream@impa.br , dream@math.haifa.ac.il , dream@tx.technion.ac.il}

\end{document}